\documentstyle{article}
\newcommand{\be}{\begin{equation}}
\newcommand{\ee}{\end{equation}}
\newcommand{\ba}{\begin{array}}
\newcommand{\ea}{\end{array}}
\newcommand{\ben}{\begin{enumerate}}
\newcommand{\een}{\end{enumerate}}

\newcommand{\ov}{\overline}
\newcommand{\tr}{{\rm tr}\, }

\newcommand{\Ldwa}{\,\stackrel{2}{\bigwedge}}

\newcommand{\w}{{\!}\wedge{\!}}
\newcommand{\Ph}{\mbox{$\bf Ph$}\, }

\font \msb=msbm10 scaled \magstep0

\newcommand{\bR}{\mbox{\msb R} }
\newcommand{\bC}{\mbox{\msb C} }

\newcommand{\ar}{\alpha }
\newcommand{\br}{\beta }
\newcommand{\gr}{\gamma }
\newcommand{\dr}{\delta }
\newcommand{\er}{\varepsilon }

\newcommand{\gd}{g^{\dagger}}

\font \eul=eufm10 scaled \magstep0
\font \eulm=eufm7 scaled \magstep0

\newcommand{\gotG}{\mbox{\eul g}}

\newcommand{\gotGm}{\mbox{\eulm g}}

\begin{document}

\title{\bf Twisted Legendre transformation}
\author{{\bf S. Zakrzewski}  \\
\small{Department of Mathematical Methods in Physics,
University of Warsaw} \\ \small{Ho\.{z}a 74, 00-682 Warsaw, Poland} }

\date{}
\maketitle

\begin{abstract}
The general framework of Legendre transformation is extended to
the case of symplectic groupoids, using an appropriate
generalization of the notion of generating function (of a
Lagrangian submanifold).
\end{abstract}

\section{Tulczyjew triple and its generalization}

The general framework of Legendre transformation was introduced
by Tulczyjew \cite{T}. It consists in recognition of the
following structure, which we call the {\em Tulczyjew triple}:
\be\label{Tt}
T^*(TQ)\; \stackrel{\ar}{\longleftarrow}\; T(T^*Q) \;
\stackrel{\br}{\longrightarrow} \; T^*(T^*Q).
\ee
Here $Q$ is the manifold of {\em configurations} (of the
system), $T^*Q$ --- its phase space (cotangent bundle),
$TQ$ --- its velocity space (tangent bundle) and $\ar ,\br$ are
natural symplectic isomorphisms. {}From those two isomorphisms,
$\br $ is the easy one: it is just the vector bundle isomorphism
induced by the symplectic form on $P=T^*Q$ ($\br $ exists for any
symplectic manifold $P$). An explicit construction of $\ar$ can
be found in \cite{T,T1}. Let us mention here how the existence of
such an isomorphism follows from the general theory of
symplectic groupoids \cite{CDW,qcp}. Cotangent bundles are
exactly those symplectic groupoids which have commutative
multiplication and connected, simply connected fibers. The
tangent bundle to such a symplectic groupoid is again a
commutative symplectic groupoid with connected, simply connected
fibers, hence a cotangent bundle to its `space of units', which in
this case coincides with $TQ$.

The dynamics is described by specifying a Lagrangian submanifold
$D\in T(T^*Q)$ (in the case of a section, we just have a
Hamiltonian vector field on $T^*Q$).  The two isomporphisms
$\ar$ and $\br$ allow to consider $D$ as a Lagrangian
submanifold of a cotangent bundle to $TQ$ and to $T^*Q$,
respectively. One can speak about generating functions of $D$ in
those two `control modes'. In the first case the function is
called the {\em Lagrangian}, and in the second case it is
(minus) the {\em Hamiltonian}. The  passage from one formulation
of dynamics to another consists in applying the {\em universal
Legendre transformation} (\ref{Tt}) to the particular case of
Lagrangian submanifold.

The aim of this paper is to stress the fact that the above
scheme can be generalized to the case when the cotangent
bundle $T^*Q$ is replaced by an arbitrary symplectic groupoid.
The set of units of the groupoid will be still denoted by $Q$.
It is now a Poisson manifold. We shall restrict our attention to
symplectic groupoids with connected and simply connected fibers
(like cotangent bundles).  They are known to be entirely
determined by their `set of units' $Q$ (with its Poisson
structure) and therefore we denote such symplectic groupoids as
$\Ph Q$. They provide a natural generalization of the phase space
(or cotangent bundle) to the case of Poisson manifolds (those which
are {\em integrable} in the sense that they are the set of units
of a symplectic groupoid). We shall call $\Ph Q$ the {\em phase
space of} $Q$, or, a {\em twisted cotangent bundle}. It is known
\cite{qcp} that any complete Poisson map of (integrable) Poisson
manifolds is lifted to a unique morphism (in the sense of
\cite{qcp}) of the corresponding twisted cotangent bundles.

We can now state the generalization of (\ref{Tt}), which may be
called the {\em twisted Tulczyjew triple}:
\be\label{tTt}
\Ph (TQ) \stackrel{\ar}{\longleftarrow} T(\Ph Q)
\stackrel{\br}{\longrightarrow} T^*(\Ph Q).
\ee
The isomorphism $\br$ is here the same as before and the
isomorphism $\ar$ arises in a similar way as before: the tangent
bundle to a symplectic groupoid over $Q$ is a symplectic
groupoid over  $TQ$. One checks easily that the resulting
Poisson structure on $TQ$ is exactly the {\em tangent} Poisson
structure \cite{C} of the original structure on $Q$. Looking at
(\ref{tTt}) we see that the dynamics $D\in T(\Ph Q)$ can be
generated --- in the usual sense --- by a Hamiltonian, which is
a function on the twisted cotangent bundle. In  Section~3
we introduce the notion of a Lagrangian submanifold of a
symplectic groupoid, {\em generated} by a function on its
Poisson manifold of units. This will allow us to consider also
the Lagrangian, which will be a function on $TQ$ (on the usual
tangent bundle, but generating in the twisted sense).

\section{Some conventions}

Recall that symplectic groupoids are symplectic-geometric models
of involutive algebras (cf. S$^*$-{\em algebras} of \cite{qcp}),
the algebra multiplication map being replaced by a symplectic
`multiplication' relation (the groupoid multiplication).  In
particular, the cotangent bundle $T^*Q$ (with its groupoid
structure) is the symplectic counterpart of the algebra of
functions on $Q$. On the other hand, from a S$^*$-algebra we get
a S$^*$-coalgebra \cite{qcp} just by inverting arrows and
passing to the dual objects (for any symplectic manifold $S$,
the role of its `dual' is played by the symplectic manifold
denoted by $\ov{S}$ which is the same manifold considered with
its opposite symplectic structure, i.e. minus the original one).
In particular, the cotangent bundle $T^*Q$ with its opposite
symplectic structure and the `co-groupoid' structure is the
symplectic counterpart of the coalgebra of measures on $Q$.
Note that we view $\Ph Q$ rather as a coalgebra than an algebra,
since in the first case the functorial correspondence between
Poisson manifolds and their phase spaces is covariant.

Now we see that it becomes important to consider on any
cotangent bundle both the canonical symplectic structure (the
exterior derivative of the canonical Liouville 1-form) and its
opposite one. As it is not clear which one should correspond to
the cotangent bundle viewed as an `algebra' (resp. `coalgebra'),
we shall adopt our first convention (it will be justified to
some extent below, when considering the phase spaces of Poisson
groups).

\vspace{1mm}

\noindent
{\bf Convention 1.} \ We shall denote by $T^*Q$ the cotangent
bundle to a manifold $Q$ with the canonical symplectic structure
and the groupoid (=`algebra') structure. By $\Ph Q$ we shall
denote the cotangent bundle with the opposite symplectic
structure and the co-groupoid (=`coalgebra') structure.

\vspace{1mm}

The next convention concerns the choice of one of the two
projections to define the correspondence between symplectic
groupoids (or rather S$^*$-coalgebras) and Poisson manifolds.

\vspace{1mm}

\noindent
{\bf Convention 2.} \ The Poisson structure on $Q$ coincides
with the image of the Poisson structure on $\Ph Q$ by the {\bf
left} groupoid projection.

\vspace{1mm}

The following convention may be viewed as a result of the two
above conventions.

\vspace{1mm}

\noindent
{\bf Convention 3.} \ If $G$ is a Lie group, then the tangent
space at the group unit, $\gotG = T_eG$, is equipped with the
Lie bracket by identifying the elements of $T_eG$ with the
corresponding fundamental vector fields of the {\em left}
translations on $G$ (i.e. the right-invariant vector fields on
$G$).

\vspace{1mm}

Indeed, the group multiplication $m\colon G\times G\to G$ is
naturally lifted to a symplectic relation from $\Ph G\times \Ph
G$ to $\Ph G$, which obeys the axioms of symplectic groupoid (it
is the {\em group S$^*$-algebra}). It defines the co-groupoid
structure on $T^*G$ (with the canonical symplectic structure).
The left groupoid projection from $T^*G$ to $\gotG$ consists in
right translation to the group unit. The Poisson
structure on $\gotG ^*$ defined by this projection is dual to
the Lie commutator on $\gotG$ coming from the right-invariant
fields.

Now we shall describe the phase spaces of Poisson groups and
motivate the first convention.

First let us recall that if $(\gotG ,\partial r)$ is a
coboundary Lie bialgebra (here $r\in \Ldwa \gotG$ and $\partial
r \colon \gotG\to\Ldwa\gotG $ is its coboundary) and $G$ is the
Lie group corresponding to $\gotG$ (convention 3), then the
Poisson (group) structure on $G$ corresponding to $\partial r$
is given by $G\ni g\mapsto rg -gr$. In this case also $g\mapsto
rg+gr$ is a Poisson structure (not a Poisson group structure
unless $r=0$).

Next, recall that for any Lie bialgebra $(\gotG,\dr )$, its {\em
Drinfeld double} is the coboundary Lie bialgebra
$(\gotG\Join\gotG ^*,\partial r_D)$, where $r_D$ is
the canonical element
\be\label{rD}
 r_D = \frac12 e_k\w e^k\in \Ldwa (\gotG\Join\gotG ^*)
\ee
($e_k$ is a basis in $\gotG$ and $e^k$ is the corresponding dual
basis in $\gotG ^*$, both considered as elements of
$\gotG\Join\gotG^*$) and $\gotG\Join\gotG^*$ is the Lie algebra
of the Manin triple (of the Lie bialgebra). The Lie group $M$
corresponding to $\gotG\Join\gotG^*$ is a Poisson group (with
the Poisson structure $M\ni\xi\mapsto r_D\xi-\xi r_D$).
Suppose that $M$ {\em decomposes} on its two subgroups $G$,
$G^*$, corresponding to Lie subalgebras $\gotG$, $\gotG ^*$,
i.e. $G\cdot G^*=M$, $G\cap G^*=\{e\}$. In this case the Poisson
group $M$ is said to be the Drinfeld double of the Poisson group
$G$ (the Poisson group structure on $G$ is defined by $\dr $).
It contains $G$ as a Poisson subgroup (here the sign of (\ref{rD})
is important!).
It is known that
the Poisson structure $\pi _+(\xi )=r_D\xi +\xi r_D$ on $M$ is
nondegenerate (everywhere on $M$) and $M$ carries a natural
structure of a (symplectic) (co-)groupoid over $G$, the left and right
projections $g$ and $g'$ of an element $\xi\in M$ on $G$ being given by
$$ \xi = g\gr = \gr ' g',\qquad \gr,\gr '\in G^* .$$
Moreover, the left projection of $\pi _+$ on $G$ is exactly the
Poisson group structure on $G$, hence $M=\Ph G$ (here also the
sign in (\ref{rD}) is important!).

For the trivial cobracket $\dr =0$, $M=T^*G$ and its symplectic
structure turns out to be opposite of the standard one, due to
the order of factors in $r_D$.

\section{Generating functions}

Let $Q$ be a manifold (with zero Poisson structure).
Any (smooth) function $f\colon Q\to \bR$ generates a Lagrangian
submanifold $L_f$ of $\Ph Q\cong T^*Q$, by the formula
\be\label{L}
L_f =\{ df(x) : x\in Q\}.
\ee
We shall give another description of this construction (its
generalization for the case of Poisson manifolds will be obvious).
Recall that the Hamiltonian vector field $X_{\phi}$ of a
function $\phi$ on a Poisson manifold $(P,\pi )$ is defined by
$$ X_{\phi}\psi = \{\phi ,\psi\} = \pi (d\phi ,d\psi ),$$
where $\psi$ is an arbitrary function on $P$. We denote by
$t\mapsto \exp (tX)$ the flow of a vector field $X$. Now observe
that $L_f$ given in (\ref{L}) equals the image of the zero section
of $T^*Q$ by $\exp (X_{f^{\rm left}})$:
\be\label{LQ}
L_f = \exp (X_{f^{\rm left}})\, Q
\ee
(here we identify the zero section with $Q$), where $f^{\rm
left}$ is the pullback of $f$ to $\Ph Q$ by the left projection
(which, of course, coincides with the right projection in the
case of the cotangent bundle). Note that here the sign of the
symplectic structure on the cotangent bundle is important:
formula (\ref{LQ}) is valid for the $\Ph Q$ symplectic form
(convention 1) equal $dx^j\w dp_j$ (using the cotangent bundle
coordinates implied by coordinates $x^j$ on $Q$).

\vspace{1mm}

\noindent
{\bf Definition.} \ Let $Q$ be an integrable Poisson manifold
and $f\colon Q\to \bR $ a complete function (i.e. a function
whose Hamiltonian vector field is complete). Then $L_f$ given by
(\ref{LQ}) is said to be the Lagrangian submanifold of $\Ph Q$
{\em generated} by $f$.

\vspace{1mm}

Note that $f^{\rm left}$ is  a complete function if and only if
$f$ is a complete function, hence the definition is correct.

\vspace{1mm}

\noindent
{\bf Lemma.} \ Formula (\ref{LQ}) is equivalent to
\be\label{Lf}
 L_f = (\Ph f)^T (\bR \times \{ 1\} ),
\ee
where $\Ph f$ is the phase lift of $f$ (i.e. the
morphism of S$^*$-coalgebras \cite{qcp} corresponding to $f$),
$(\Ph f)^T$ denotes the transposed relation of $\Ph f$ and
$$\bR \times \{ 1\} = \{ (e,t)\in \Ph \bR : e\in \bR, t=1\}$$
is the Lagrangian submanifold of $\Ph \bR$ given by the `time'
$t=1$.

\vspace{1mm}

\noindent
{\em Proof.} \ \ Let us calculate $\Ph f$ by the method of
characteristics \cite{qcp}. We consider the coisotropic submanifold
in $\Ph Q\times \overline{\Ph \bR}$ given by the single constraint
$f^{\rm left} (\xi) = e$. The characteristics are the integral
curves of
$$ X_{f^{\rm left} - e}=X_{f^{\rm left}} + X_e = X_{f^{\rm
left}} +\frac{\partial}{\partial t}.$$
An initial condition of the form $(x;f(x),0)$
evolves like $t\mapsto (\exp (tX_{f^{\rm left}})x;f(x),t)$, hence
the graph of $\Ph f$ is
$$ \{ \exp (tX_{f^{\rm left}})x; f(x),t) : x\in Q, t\in \bR \}.$$
This implies (\ref{Lf}).

\hfill $\Box$

Since $\Ph f$ can be calculated also using the right
projections, from (\ref{Lf}) it follows that in (\ref{LQ}) we
can replace $f^{\rm left}$ by $f^{\rm right}$ (the pullback of
$f$ to $\Ph Q$ by the right projection) and the result will be
the same.

Since $X_{f^{\rm left}}$ is right-invariant \cite{CDW}, we have
the following property:
\be\label{Lfxi}
\exp (X_{f^{\rm left}})\xi = L_f\cdot \xi\qquad \mbox{for}\;\;
\xi\in \Ph Q
\ee
(the dot on the right hand side is the groupoid multiplication).
Similarly,
\be\label{Lfxi2}
\exp (X_{f^{\rm right}})\xi =\xi\cdot L_f \qquad \mbox{for}\;\;
\xi\in \Ph Q .
\ee

\vspace{1mm}

\noindent
{\bf Example 1.} \  Let $\gotG$ be a Lie algebra.
 Any $X\in \gotG$ is a linear function on $Q=\gotG ^*$. It generates
a Lagrangian submanifold in $\Ph \gotG^* = T^*G$, where $G$ is
the (connected and simply connected) Lie group corresponding to
$\gotG $. We have
$$ f^{\rm left}(\xi ) = \left\langle X^{\rm right\; invariant},\xi
\right\rangle =  \left\langle X_{\rm left \; translation},\xi
\right\rangle ,$$
therefore $f^{\rm left}$ is just the canonical moment map of the
left translations, evaluated on $X$. It follows that $f^{\rm
left}$ generates left translations  and $L_f=\exp (X_{f^{\rm
left}}) T_e^*G = T^*_{\exp X}G$ is the fiber of the cotangent
bundle at $g=\exp X$.

\vspace{1mm}

\vspace{1mm}

\noindent
{\bf Example 2.} \
Let $S$ be a simply connected symplectic manifold. Its phase
space is the {\em pair} co-groupoid $\Ph S = S\times \ov{S}$
with $S\cong \{ (x,x) : x\in S\}$ as the set of units.
The left and right projections are
$$ (x,y)\mapsto x \qquad\mbox{and}\qquad (x,y)\mapsto y,$$
respectively, hence for any function $f$ on $S$ we have
$$ f^{\rm left} (x,y)=f(x),\qquad f^{\rm right}(x,y)=f(y),
$$
and
$$ X_{f^{\rm left}} = X_f \oplus 0,\qquad X_{f^{\rm right}} =
0 \oplus (-X_f).
$$
Assuming that $f$ is complete we have then
$$ \exp (X_{f^{\rm left}}) (x,x)= (\exp (X_f)x,x),\qquad
\exp (X_{f^{\rm right}}) (x,x)= (x,\exp (-X_f)x),$$
and it follows that the Lagrangian submanifold generated by $f$
in $\Ph S$ is the graph  of the symplectic transformation
$\exp X_f \colon S\to S$.

\vspace{1mm}

The last example shows that in general it is very difficult to
find the generating function of a given Lagrangian submanifold:
one has to find a Hamiltonian flow which includes a given
symplectic transformation (the generator needs not be unique;
it is seen also in Example 1).
Also the generation procedure is in general quite ineffective,
because it requires to solve the equations of motion for some
Hamiltonian. One can hope it becomes simpler for the Casimir
functions of a given Poisson structure (functions which Poisson
commute with all other functions, i.e. functions constant on
symplectic leaves). Let us discuss shortly this case. Note that each Casimir
function is complete. For such functions $f^{\rm
left}=f^{\rm right}$ and $X_{f^{\rm left}}$ is tangent both to
left and right fibers. Let $Q_0\subset Q$ be the (open
and dense) subset of points of the maximal rank of the Poisson
structure. The corresponding part $\Ph Q_0$ of $\Ph Q$
carries a natural coisotropic foliation (level sets of the
pullbacks of Casimir functions). The leaves of the
corresponding characteristic isotropic foliation are spanned by
the action of the Hamiltonian vector fields of the pullbacks of
Casimir functions. Each such characteristic leaf must belong to the
intersection of a left fiber and a right fiber. In fact the
intersection coincides with the characteristic leaf, because the
intersection of tangent spaces to left and right fibers are exactly
the tangent spaces to the characteristic foliation.
In particular, the intersection of the left and the right fiber over
the same point $x\in Q_0$ --- so called {\em isotropy group} of
$x$ (cf. \cite{CDW}) --- is spanned by the action on $x$ of flows of
the  Hamiltonian vector fields of pullbacks of Casimir
functions. Let $\xi $, $\eta$ be two elements of the isotropy
group which can be written as
$$ \xi = \exp (X_{f^{\rm left}}) x= L_f\cdot x,\qquad
\eta = \exp (X_{g^{\rm left}}) x= L_g\cdot x, $$
for some Casimir functions $f,g$ (formula (\ref{Lfxi})). We have
\be\label{comm}
\xi \cdot \eta = x\cdot L_f\cdot L_g\cdot x =
x\cdot L_g\cdot L_f \cdot x = \eta\cdot \xi
\ee
because $L_f\cdot L_g =L_g\cdot L_f$ for Casimir functions. The
latter property follows from
$$ L_f\cdot L_g\cdot\xi = \exp (X_{f^{\rm left}})\exp (X_{g^{\rm
right}})\xi = \exp (X_{g^{\rm right}}) \exp (X_{f^{\rm
left}})\xi = L_g\cdot L_f\cdot\xi .$$
{}From (\ref{comm}) it follows that the isotropy group is abelian.
Since $X_{f^{\rm left}}$ is right invariant, formula
$$(X_{f^{\rm left}})|_x\mapsto  \xi = \exp (X_{f^{\rm left}}) x $$
actually defines the exponential map for this group. We conclude
that in order to calculate $L_f$, it is sufficient to calculate
$df(x)$ at each point $x\in Q_0$, and use the groupoid structure
restricted to the isotropy subgroup of $x$:
$$ L_f|_{Q_0}=\{ \exp ( (df (x))^{\sharp }) x :  x\in Q_0\}
$$
(here $\sharp$ denotes  `raising of indices' by the Poisson
structure on $\Ph Q$).

We remark that our result on the commutativity of the `minimal'
isotropy subgroups (isotropy subgroups of maximal symplectic
leaves)  in the case of the `group S$^*$-algebra' (cf. Section~2)
means that isotropy subgroups of elements belonging to maximal
coadjoint orbits are abelian (cf. \cite{X}). In the case of
`Poisson group S$^*$-algebra' (i.e. the symplectic groupoid
$G\cdot G^*$ over $G^*$), it means the same, with coadjoint
action replaced by the dressing one.

\vspace{1mm}

\noindent
{\bf Example 3.} \ Let $V$ be a vector space with a constant
Poisson structure $r$, i.e. the Poisson brackets of linear
coordinates $x^j$ are
$$ \{ x^j ,x^k\} = r^{jk}= - r^{kj},$$
where $r^{jk}$ are some constants (the coefficients of $r$ in
the corresponding basis). The phase space $\Ph V$ can be
realized \cite{abel} in $T^*V=V\oplus V^*$ with the canonical
structure
$$ \{ q^j,q^k\}=0,\qquad\{ q^j,p_k\} = \dr ^j_k,\qquad
\{ p_j,p_k\} =0$$
and the left projection equal
$$x = (q,p)_L= q +\frac12 r (p), \qquad [r(p)]^j := p_kr^{kj} .$$
For any {\em linear\/} function $f(x) = a_jx^j$, we can find the
flow of $f^{\rm left}(q,p)=a_j (q^j +\frac12 p_kr^{kj}$). We have
$$ \dot{q}^k = \{ f^{\rm left} ,q^k\} = \frac12 a_jr^{jk},\qquad
\dot{p}_j = a_j ,$$
i.e. $x^k(t)=x^k (0) +\frac12 a_jr^{jk}t$, $p_j = p_j (0)+a_jt$,
hence
$$ L_f = \{ (q,p) : p=a\} $$
is the same `horizontal' Lagrangian plane as in the usual case.

If $f$ is {\em quadratic}, then $f^{\rm left}$ is also quadratic
and $\exp (X_{f^{\rm left}}) $ is a linear symplectic
transformation, hence $L_f$ is a Lagrangian (linear) subspace of
$V\oplus V^*$, as in the usual case.

\section{Particle on the Poisson $SU(N)$}

An explicit calculation of the relation between the Lagrangian
and the Hamiltonian may be very difficult in the `twisted' case,
due to the nontrivial generation procedure.
In this section we analyze this problem a little bit for
an analogue of the free motion, when the Poisson configuration
manifold $Q$ is the standard Poisson $SU(N)$ group.
The description of $\Ph Q$ is given in \cite{sun}: it is
$SL(N,\bC )$ --- the Manin group (with the $\pi _+$ Poisson
structure), which contains both $G=SU(N)$ and its (Poisson) dual $G^*=SB(N)$.
We can ask two questions:
\ben
\item Which Lagrangian corresponds to the Hamiltonian
\be\label{ham}
 {\cal H} (g) : = \frac12 \tr g^{\dagger}g,\qquad g\in SL(N,\bC),
\ee
introduced in \cite{sun} and describing a natural analogue of
the free motion on $SU(N)$?
\item Which Hamiltonian corresponds to the `non-deformed free' Lagrangian
\be\label{lag}
 {\cal L} (v) := \frac12 v^2,\qquad  v\in TQ,
\ee
where the square of $v\in TQ=TG$ is taken with respect of a bi-invariant
metric on $Q=G$?
\een

Unfortunately, we are still premature to answer these questions.
In the sequel we show how far we are able to reduce the problem.

The dynamical equations implied by (\ref{ham}) were obtained in
\cite{sun} and look as follows
\be\label{eqmot}
 \dot{g} = i\er [g\gd g - \frac1{N} (\tr \gd g) g].
\ee
If we decompose $g=u\gr $ with $u\in G$, $\gr \in G^*$, we get
\be\label{dyn}
 \dot{\gr}=0,\qquad u^{-1}\dot{u}= F(\gr ):= i\er [\gr\gr ^{\dagger}
 - \frac1{N} (\tr \gr \gr ^{\dagger})].
\ee
It was shown in \cite{sun} that $F$ is bijective, hence
the dynamics (\ref{eqmot}) is a section of the left projection
$$ T\, SL(N,\bC)\ni \dot{g} =\dot{u}\gr +u\dot{\gr}\mapsto
\dot{u} \in TG,
$$
namely
$$
TG\ni \dot{u}\mapsto \dot{g}=\dot{u}F^{-1} (u^{-1}\dot{u})\in
T\, SL(N,\bC ).
$$
Similarly, decomposing $g$ in the opposite order, $g=\gr u$, we get
\be\label{dynr}
 \dot{\gr}=0,\qquad \dot{u}u^{-1} = E (\gr ):= i\er [\gr ^{\dagger}\gr
 - \frac1{N} (\tr  \gr ^{\dagger}\gr )].
\ee
Since $E$ is also bijective,
the dynamics (\ref{eqmot}) is also a section of the right projection
$$ T\, SL(N,\bC)\ni \dot{g} =\gr\dot{u} +\dot{\gr}u\mapsto
\dot{u} \in TG,
$$
namely
\be\label{hamdyn}
TG\ni \dot{u}\mapsto \dot{g}=E ^{-1} (\dot{u} u^{-1})\dot{u} \in
T\, SL(N,\bC ).
\ee
In order to find its generating function on $TQ$  --- the
Lagrangian --- we have to try to calculate the Lagrangian
submanifolds generated by some `natural' (class of) Lagrangian(s)
and see whether this fits the original dynamics. Thus we are in
fact led to the Question 2 above.

In order to find the Lagrangian submanifold generated by
a function on $TG$, we should have a formula for the left
projection from $T \Ph Q$ to $TQ$ and an effective control of
the symplectic structure of $T\Ph Q$. For this latter purpose,
it is more convenient to work in the canonical symplectic
structure of $T^* \Ph Q$. We can use it passing from $T\Ph Q$ to
$T^* \Ph Q$ with the help of $\br$ (formula (\ref{tTt})).
Fortunately, the isomorphism $\br $ not only preserves the
symplectic structure, but also transforms the tangent groupoid
into cotangent groupoid of $\Ph Q$. Indeed, as can be easily seen,
for any symplectic relation $\rho\colon X\to Y$ (from a
symplectic manifold $X$ to a symplectic manifold $Y$), $T\rho$
is transformed into $\Ph \rho$ by $\br$. We just apply this fact
to the case of the groupoid multiplication relation in $\Ph Q$.
It follows that the left projection in $T\Ph Q$ over $TQ$ is
just the left projection in $T^*\Ph Q$ over $(TQ)^{\circ}$ (the
conormal bundle to $TQ$ in $T^*\Ph Q$), composed with the
inverse of $\br$.

The left projection in $T^* \Ph Q$ over $(TQ)^{\circ}$ for $Q$
being a double Lie group has been calculated in \cite{qcp}
(formula (25)). It is given by
\be\label{lp}
 \xi _L =[P(\xi (u\gr ) ^{-1})]u \qquad
\mbox{for}\;\; \xi\in T^*_{u\gr }(G\cdot G^*),
\ee
where $P$ is the projection on $\gotG ^{\circ}$  parallel to
$(\gotG ^*)^{\circ} $ (here $u\in G$, $\gr \in G^*$).
Let $s$ be the canonical symmetric map from $(\gotG \Join
\gotG^* )^* $ to $\gotG\Join \gotG ^*$. We have
$$ [(su) \xi _L]u^{-1} = s (\xi _L u^{-1}) =
s(P(\xi g^{-1}))= \xi g^{-1}|_{\gotGm ^*}\in \gotG
\qquad \mbox{for}\;\; g=u\gr\in G\cdot G^*.$$
Using the formula for $\sharp =\br ^{-1}$,
$$ \sharp = \frac12 (Rg+gR) (gs)= \frac12 (Rg+gR)(sg)$$
(which follows from $\pi _+=r_D g+gr_D$), where $R$ is the
reflection in $\gotG ^*$ parallel to $\gotG$, we obtain
$$
\sharp \xi _L = \frac12 (Ru+uR) (su)\xi _L =
- (su)\xi _L = -[\xi g^{-1}|_{\gotGm ^*} ]u.
$$
If ${\cal L}\colon TG\to \bR$ is right-invariant, then
$$ {\cal L}(\sharp \xi_L)=
{\cal L} (-(\xi g^{-1})|_{\gotGm ^*})$$
is a right-invariant function on $T^*(G\cdot G^*)$. For
instance, if ${\cal L}$ is given by a right-invariant metric  on $G$,
such that for  $v\in \gotG$ we have
$${\cal L} (v) = \frac12 \left\langle l,v\otimes v\right\rangle ,$$
where $l$ is a symmetric element of $\gotG ^*\otimes \gotG ^*$, then
$${\cal L} (\sharp \xi _L) = \frac12 \left\langle l,
\xi g^{-1}|_{\gotGm ^*}\otimes \xi g^{-1}|_{\gotGm ^*}\right\rangle
=\frac12 \left\langle l,\xi g^{-1}\otimes \xi g^{-1}\right\rangle
=\frac12 \left\langle lg,\xi \otimes \xi \right\rangle ,
$$
so we are interested in the flow of the right-invariant
contravariant metric on $G\cdot G^*$ which at the unit is
tangent to $\gotG ^*$ (and equal $l$). Of course, in such a
situation, the left translation of $\xi $ to the group unit is constant:
$$ g^{-1}\xi = \xi _0\in (\gotG\Join\gotG^*)^*.$$
Since we are considering the flow of a function which is
constant on the fibers of the left projection, we know that the
right projection, given by the formula analogous to (\ref{lp}),
\be\label{rp}
 \xi _R =uP((\gr u)^{-1}\xi ) \qquad
\mbox{for}\;\; \xi\in T^*_{\gr u}(G\cdot G^*),
\ee
is preserved. It follows that $u$ is constant and
$$ \xi = g\xi _0 = \gr (t) u \xi _0.$$
We are interested in the initial conditions $\gr (0)=e$, $\xi
_0\in \gotG ^{\circ}$, hence
\be\label{solu}
 \xi = \gr (t) \eta _0 u
\ee
for some $\eta _0\in \gotG ^{\circ}$. Due to the right
invariance, it is sufficient to solve the flow for the initial
condition being just $\eta _0\in \gotG^{\circ}\cong \gotG \cong
(\gotG ^*)^*$ (i.e. $u=e$), hence $\xi = \gr (t)\eta _0$. Now it
is easy to see that the problem is reduced to solving the flow
on
$$ T^*G^*\cong G^*\cdot \gotG ^{\circ}\subset T^*(G\cdot G^*),$$
generated by the right-invariant contravariant metric on $G^*$,
defined by $l$, for the initial condition $\eta _0\in \gotG $.
Note that the flow of a right-invariant positive contravariant
metric is always complete on the cotangent bundle, because such
a Hamiltonian is a pullback (by the left projection) of a
function on the dual of the Lie algebra which has compact level
sets (hence the function is complete).

We conclude that the whole information
about the  Lagrangian submanifold generated by (\ref{lag}) is 
contained in the `one-fiber' Legendre transformation
$$ \gotG \cong \gotG ^{\circ}\ni \eta _0 \mapsto \Phi (\eta
_0)\in G^* ,$$
where $\Phi (\eta _0) =\gr (1)$ is the `time equal 1' point of
the geodesic $t\mapsto \gr (t)$
which starts at $t=0$ with velocity corresponding to $\eta _0$.
The Lagrangian submanifold generated by ${\cal L}$ is
described by
\be\label{Lagdyn}
\xi = \Phi (\eta _0)\eta _0 u,\qquad \eta _0\in \gotG
^{\circ}\cong \gotG,\; u\in G.
\ee
In order to compare it with (\ref{hamdyn}), let us note that
$\dot{g}\in g\gotG \cap \gotG g$ (because both left and right
projections of $\dot{g}$ on $G^*$ are zero), therefore
$$ \frac12 (Rg+gR) \dot{g} = - \dot{g} .$$
It follows that
$$\br (\dot{g}) = - (sg)(\dot{g})= -(sg) [E^{-1} (\dot{u}u^{-1})\dot{u}]
$$
$$ = -(sg) [E^{-1} (v)vu] = - E^{-1} (v) s(v) u ,$$
where $v := \dot{u}u^{-1}\in \gotG $, hence we see that ${\cal
L}$ and ${\cal H}$ generate the same dynamics when
\be\label{redu}
 -E^{-1} = \Phi .
\ee
Thus we have reduced the problem of comparing the Hamiltonian
and the Lagrangian dynamics to the equality of the above type.
Some straightforward calculations for $SU(2)$ show that the
equality does not hold for ${\cal H}$ and ${\cal L}$ as given by
(\ref{ham}) and (\ref{lag}), respectively. The above reduction
of the problem to the form (\ref{redu}) is valid however for
${\cal H}$ coming from any Casimir on $G^*$ and ${\cal L}$
coming from any invariant on $\gotG$. Therefore there is still
room for finding appropriate compatible pairs $({\cal H},{\cal
L})$.

\end{document}